\documentclass[10pt]{article}

\usepackage{theorem,amssymb,amsmath}

\usepackage{graphicx}

\usepackage{color}                    
\definecolor{red}{rgb}{1,0,.2}        
\definecolor{cjp}{rgb}{.1,.7,.2}        
\definecolor{fmdc}{rgb}{1,0,.8}        
\definecolor{gree}{rgb}{0,0.7,0.2}        
\newcommand*{\gr}[1]{{\color{gree} #1}} 
\definecolor{blue}{rgb}{0,0,.9}

\newcommand*{\bl}[1]{{\color{blue} #1}} 

\newcommand{\sF}{\sigma_{\rm \scriptstyle F}}
\newcommand{\xG}{x_g}

\newcommand{\B}{{\rm B}}

\newcommand{\E}{{\rm E}}

\newcommand{\bC}[1]{\bl{$\bigcirc$\!\!\!\!{\footnotesize $#1$}}}

\newcommand{\pq}[2]{ \left( \begin{smallmatrix} #1\\ #2 \end{smallmatrix} \right)}                     
\newcommand{\PQ}[2]{ \left( \begin{matrix} #1\\ #2 \end{matrix} \right)}                               

\def \endpf{{\ \ $\Box$ \medbreak}}

\newtheorem{theorem}{Theorem}

\newtheorem{lemma}[theorem]{Lemma}

\newtheorem{proposition}[theorem]{Proposition}
\theorembodyfont{\rm}

\advance \topmargin by -\headheight
\advance \topmargin by -\headsep
\evensidemargin \oddsidemargin
\begin{document}

\begin{centering}

{\Large \textbf{The sum of digits functions of the Zeckendorf and the base phi expansions}}


\bigskip

{\bf \large Michel Dekking}

\medskip

{ DIAM,  Delft University of Technology, Faculty EEMCS,\\ P.O.~Box 5031, 2600 GA Delft, The Netherlands.}

\medskip

{\footnotesize \it Email:  F.M.Dekking@TUDelft.nl}

\end{centering}

\smallskip


\smallskip

\begin{abstract}
 \noindent We consider the sum of digits functions for both base phi, and for the Zeckendorf expansion of the natural numbers. For both sum of digits functions we present morphisms on infinite alphabets such that these functions viewed as infinite words are letter-to-letter projections of fixed points of these morphisms. We characterize the first differences  of both  functions  a) with generalized Beatty sequences, or unions of generalized Beatty sequences, and b) with morphic sequences.
\end{abstract}

\medskip


\vspace*{2cm}



\section{Introduction}

Our  interest in this paper is in the irrational base phi, where  phi is the golden mean $\varphi=(1+\sqrt{5})/2$, and in the Zeckendorf representation of the natural numbers.

For both expansions we give in Theorem \ref{th:Zeck-inf}, respectively Theorem \ref{th:phi-inf} a morphism on an infinite alphabet, such that the sum of digits functions of these expansions considered as infinite words are letter-to-letter projections of fixed points of these morphisms.

We then will show how these results permit to give precise information on the  first differences of the sum of digits functions.
The {\it first differences} of a function $f:\mathbb{N}_0\rightarrow\mathbb{N}_0$\; are given by the function $\Delta f$ defined by $$\Delta f(N)=f(N+1)-f(N),\quad \text{for}\; N=0,1,2\dots.$$
    We shall focus on the signs of $\Delta f$.
   A number $N$ is called a {\it point of increase} of a function $f:\mathbb{N}_0\rightarrow\mathbb{N}_0$\; if \;$\Delta f(N)>0$. It is called a {\it point of constancy} if \;$\Delta f(N)=0$, and a {\it point of decrease} if \;$\Delta f(N)<0$.

   For base two it is simple to see that  the points of increase of the sum of digits function $s_2$ with values $0,1,1,2,1,2,2,3,\dots$ are given by the even numbers, and that the points of constancy and decrease are given by the numbers $1\! \mod 4$, respectively $3\! \mod 4$.

\medskip

We will prove (Theorem \ref{th:Zeck-GBS} and Theorem \ref{th:phi-GBS}) for both the Zeckendorf representation and the base  phi expansion that the points of increase, constancy and decrease are all given by unions of generalized Beatty sequences, as studied in \cite{GBS}.
These are  sequences  $V$ of the type  $$V(n) = p\lfloor n \alpha \rfloor + q\,n +r , \quad n\ge 1,$$
 where $\alpha$ is a real number, and $p,q,$ and $r$ are integers. Here we denoted the floor function by $\lfloor \cdot \rfloor$.

We will also prove that the  first differences of the sequences of points of increase, constancy and decrease are all morphic sequences. See Theorem \ref{th:Zeck-morph}  for the Zeckendorf representation, and  Theorem  \ref{th:phi-morph} for the base phi expansion.

\medskip

  A prominent role in this paper, both for base phi and the Zeckendorf expansion, is played by $(\lfloor n\varphi\rfloor)$,  the well known lower Wythoff sequence.

A standard result (see, e.g., \cite{lothaire}) is   that the sequence $\Delta(\lfloor n\varphi\rfloor)$  is equal to the Fibonacci word $x_{1,2} = 1211212112\dots$  on the alphabet $\{1,2\}$, i.e., the unique fixed point of the morphism $1\rightarrow 12, 2\rightarrow 1$. More generally, we have  the following simple lemma.

\begin{lemma}\label{lem:diff}{ \rm \bf(\cite{GBS})} Let $V = (V(n))_{n \geq 1}$ be the generalized Beatty
sequence defined by $V(n) = p\lfloor n \varphi \rfloor + q n +r$, and let $\Delta V$ be the
sequence of its first differences. Then $\Delta V$ is the Fibonacci word on the alphabet
$\{2p+q, p+q\}$. Conversely, if $x_{a,b}$ is the Fibonacci word on the alphabet
$\{a,b\}$,  then any $V$ with $\Delta V= x_{a,b}$ is a generalized Beatty sequence
$V=((a-b) \lfloor n \varphi \rfloor)+(2b-a)n+r)$ for some integer $r$.
\end{lemma}

Let $A(n)=\lfloor n \varphi \rfloor$, and $B(n)=\lfloor n \varphi^2 \rfloor$. It is well known that $A$ and $B$ form a pair of Beatty sequences, i.e., they are disjoint with union $\mathbb{N}$. In the next lemma, $V\!A$ is the composition given by $V\!A(n)=V(A(n))$.

\begin{lemma}\label{lem:split}{ \rm \bf(\cite{GBS})} Let $V$ be a generalized Beatty sequence given by $V(n) = p\lfloor n \varphi \rfloor + q n +r $, $n\ge 1$. Then $V\!A$ and $V\!B$ are generalized Beatty sequences with parameters $(p_{V\!A},q_{V\!A},r_{V\!A})=(p+q,p,r-p)$ and $(p_{V\!B},q_{V\!B},r_{V\!B})=(2p+q,p+q,r)$.
\end{lemma}

\section{The Zeckendorf sum of digits function }\label{sec:Zeck}

  Let $F_0=0, F_1=1, F_2=1,\dots$ be the Fibonacci numbers.
  Ignoring leading and trailing zeros, any  natural number $N$ can be written uniquely   with digits $d_i=0$ or 1, as
  $$N= \sum_{i\ge0} d_i F_{i+2},$$
  where $d_id_{i+1} = 11$ is not allowed.     We denote the Zeckendorf expansion of $N$ as $Z(N)$, with digits $d_i(N)$.

  Let $s_Z$ be the sum of digits of such an expansion: for $N\ge 0$
   $$s_Z(N)=\sum_{i\ge0} d_i(N).$$
 We have
  $$(s_Z(N))=(0,1,1,1,2,1,2,2,1,2,2,2,3,1,2,2,2,3,2,3,3,1,2,2,2,\dots)$$

\medskip

Our first result is that $s_Z$ is a morphic sequence.

\smallskip

\begin{theorem}\label{th:Zeck-inf}
  The function $s_Z$, as a sequence, is a morphic sequence on an infinite alphabet, i.e., $(s_Z(N))$ is a letter to letter projection of a fixed point of a morphism $\tau$.
The alphabet is $\{0,1,...,j,...\}\times\{0,1\}$, and $\tau$ is the morphism given by\\[-.6cm]
\begin{align*}
      \tau\Big(\PQ{j}{0}\Big)& \;=\; \PQ{j}{0}\, \PQ{j\!+\!1}{1}, \\
      \tau\Big(\PQ{j}{1}\Big)& \;=\; \PQ{j}{0}.
\end{align*}
The letter-to-letter map is given by the projection on the first coordinate: $\pq{j}{i}\rightarrow j$ for $i=0,1.$
The fixed point $x_\tau$ of $\tau$ with initial symbol $\pq{0}{0}$ projected on the first coordinate equals $(s_Z(N))$.
\end{theorem}

\medskip

\noindent {\it Proof:} See the Comments of sequence A007895 in \cite{oeis} for a proof of this. \hfill\endpf

\medskip


\medskip

Let $I_Z, C_Z$ and $D_Z$ be the functions listing the points of increase, constancy, and decrease of the function $s_Z$.
We have
$$I_Z=(0,3,5,8,11,13,16,\dots)\footnote{$I_Z$ is the sequence  A026274 in \cite{oeis}.},\; C_Z=(1,2,6,9,10,14,\dots),\; D_Z=(4,7,12,17,20,25,\dots).$$

To state our results it is actually convenient to define $D_Z=(-1,4,7,12,17,20,25,\dots)$.

\medskip

When $(a_n)$ and $(b_n)$ are two increasing sequences, indexed by $\mathbb{N}$, then we mean by the union of $(a_n)$ and $(b_n)$ the increasing sequence whose terms go through the set $\{a_n, b_n: n\in \mathbb{N}\}$.

\medskip

\begin{theorem}\label{th:Zeck-GBS} The function $I_Z$, the points of increase of the function $s_Z$,  is given for $n=1,2,\dots$ by
$$I_Z(n)=\lfloor n\varphi \rfloor + n-2.$$
The function $C_Z$, the points of constancy of the function $s_Z$,  is given for $n=1,2,\dots$ by the union of the two generalized Beatty sequences with terms
$$2\lfloor n\varphi \rfloor + n-2\quad\text{and}\quad 3\lfloor n\varphi \rfloor + 2n-3.$$
The function $D_Z$, the points of decrease of the function $s_Z$,  is given for $n=1,2\dots$ by
$$D_Z(n)=2\lfloor n\varphi \rfloor + n-4.$$
 \end{theorem}

 \medskip

 \noindent {\it Proof:}
Let $I_Z$ be the sequence of the points of increase of the function $s_Z$.

Projection on the second  coordinate of $\tau$ yields the Fibonacci morphism $\sF$ given by
$$\sF(0) = 01,\; \sF(1)=0.$$
Thus the second coordinates of the  fixed point of $\tau$ equal the infinite Fibonacci word $x_{0,1}$ = 0100101001001....
Obviously, the increase points of $s_Z$ occur if and only if the word $(j,0)\, (j\!+\!1,1)$ occurs in the fixed point $x_\tau$ of $\tau$ if and only if the word 01 occurs in $x_{0,1}$. Since 11 does not occur in $x_{0,1}$, this means that we have to shift the positions of 1's in $x_{0,1}$ by 1. It is well known that the positions of 1  are given by the upper Wythoff sequence $(\lfloor n\varphi^2 \rfloor )=(\lfloor n\varphi \rfloor + n)$. Since the first coordinate of the fixed point of $\tau$ starts from index 0, and the second from index 1, we have to replace $n$ by $n+1$, and this yields the first result of Theorem \ref{th:Zeck-GBS}.

\medskip

The points of constancy are more difficult to characterize with the fixed point $x_\tau$ than the points of increase.
We therefore take another approach. Write $Z(N)=\dots w$, where $w$ is a word of length 4. Then $w$ can be any word of the $0$-$1$-words of length 4 containing no 11. Obviously, the three words $w=0000, w=0100$ and $w=1000$  give points of increase.

Furthermore the numbers $N$ with $Z(N)$ ending in $w=0001, 1001$ and $w=0010$ give
 $$Z(N)=\dots001 \Rightarrow Z(N+1)\doteq \dots002\doteq\dots010,\quad Z(N)=\dots0010 \Rightarrow  Z(N+1)\doteq \dots0011\doteq\dots0100.$$
We see that these give  points of constancy.

Finally, we show that the $N$ with $Z(N)$ having suffix $w=0101$ or $w=1010$ give points of decrease.
\begin{align*}
Z(N)&=\dots0101\;\;\quad \Rightarrow \quad Z(N+1)\doteq \dots0102\doteq\dots0110\doteq\dots1000,\\
Z(N)&=\dots01010  \quad \Rightarrow \quad Z(N+1)\doteq \dots01011\doteq\dots01100\doteq\dots10000.
\end{align*}
In both cases at least one digit 1 is lost, so these $N$ are the points of decrease.

With this knowledge we can apply  Theorem 2.3 and Proposition 2.8 in the paper \cite{Dekking-Zeck},  obtaining that one part of $I_Z$ is given by the generalized Beatty sequence
$(2\lfloor n\varphi\rfloor+n-2)$ and the other part is given by $(3\lfloor n\varphi\rfloor+2n-3)$.

 Again from Theorem 2.3 and Proposition 2.8 in the paper \cite{Dekking-Zeck}, we obtain that $(D_Z(n+1))$ is the union of the two generalized Beatty sequences $(3\lfloor n\varphi\rfloor+2n-1)$ and $(5\lfloor n\varphi\rfloor+3n-1)$.

  It is not a simple matter to see that this union is given by the single generalized Beatty sequence $(2\lfloor n\varphi \rfloor + n-4)$, where the index starts at $n=2$.

   Let us write $V(p,q,r)=(p\lfloor n\varphi \rfloor + qn+r)_{n\ge 1}.$
We have proved so far that $I_Z=V(1,1,-2)$, and $C_Z$ is the union of $V(2,1,-2)$ and  $V(3,2,-3)$. If we add  2 to all terms of these sequences, we obtain the three sequences $V(1,1,0)$,  $V(2,1,0)$, and  $V(3,2,-1)$.

The triple of sequences
$$\{V(1,1,0),  V(2,1,0), V(1,1,-1)\}$$
 is known as the `first classical complementary triple', i.e., these are three disjoint sequences with union $\mathbb{N}$. See page 334 in \cite{GBS}. The third sequence of this triple, $V(1,1,-1)$, can be written as a disjoint union of the two sequences $V(3,2,-1)$ and $V(2,1,-2)$, by Lemma \ref{lem:split}. Thus

$$\{V(1,1,0),  V(2,1,0), V(3,2,-1), V(2,1,-2)\}$$
forms a complementary quadruple. If we subtract 2 from all terms of these four sequences, the first gives $I_Z$, the second and the third together, $C_Z$. Since $\{I_Z,C_Z,D_Z\}$ is a complementary triple, with union $\{-1,0,1,2,\dots\}$ this implies that $(D_Z(n+1))$  has to be equal to $V(2,1,-4)$. \hfill\endpf

\medskip

 Next, we give a characterization of $I_Z, C_Z$ and $D_Z$ in terms of morphisms.

 \begin{theorem}\label{th:Zeck-morph} The points of increase of the function $s_Z$  are given by the sequence $I_Z$, which has $I_Z(1)=0$, and $\Delta I_Z$  is the fixed point of the Fibonacci morphism $3\rightarrow 32,\; 2\rightarrow 3$.

The points of constancy of the function $s_Z$  are given by the sequence $C_Z$, which has $C_Z(1)=1$, and $\Delta C_Z$  is the fixed point of the 2-block Fibonacci morphism on the alphabet $\{1,4,3\}$ given by  $1\rightarrow 14,\; 3\rightarrow 14,\; 4\rightarrow 3$.

The points of decrease of the function $s_Z$  are given by the sequence $D_Z$, which has $D_Z(1)=-1$, and $\Delta D_Z$  is the fixed point of the Fibonacci morphism $5\rightarrow 53,\; 3\rightarrow 5$.
\end{theorem}

\medskip

For the proof of Theorem \ref{th:Zeck-morph} we have to make some preparations. Let $\Lambda_3:=\{2\}, \Psi_3:=\{0,1\}=:[0,1]$, and define for $n\ge 4$ the intervals of integers $\Lambda_n$ and $\Psi_n$ by
$$\Lambda_n:=[F_n,F_{n+1}-1], \; \Psi_n:=[0,F_{n}-1].$$
The $(\Lambda_n)$ form a partition of $\mathbb{N}_0\setminus \{0,1\}$, and the $(\Psi_n)$ satisfy
\begin{equation}\label{eq:rec}
\Psi_{n+1}=\Psi_n\cup \Lambda_n.
\end{equation}
For an interval $I$, let $C_Z(I)$ denote the points of increase lying in the interval $I$.  Also, let $\Delta C_Z(I)$ denote the first differences of the points of increase lying in the interval $I$, considered as a word on the alphabet $\{1,2,3,4\}$. At first sight, the latter definition is problematic, as one has to know the first point of increase after the last element of  $C_Z(I)$. However, we shall only consider intervals $I=\Lambda_n$ and $I=\Psi_{n+1}$, which both are followed by $\Lambda_{n+1}$, and one verifies easily that the first point of increase in $\Lambda_{n+1}$ is always the second point. Actually, this follows directly from the following lemma.
\medskip

\begin{lemma}\label{lem:cZ} For all $n\ge 3$ one has $C_Z(\Lambda_{n+1})=C_Z(\Psi_n)+F_{n+1}$. \end{lemma}

\noindent {\it Proof:} We used the notation $A+y=\{x+y:x\in A\}$ for a set $A$, and a number $y$.
 The lemma follows from the basic Zeckendorf recursion: the numbers $N$ in $\Lambda_{n+1}$ all have a digit 1 added to the expansion of $N-F_{n+1}$. \hfill\endpf

 \medskip

 Let $h$ be the morphism on the alphabet $\{1,3,4\}$ given by
$$h(1)= 14,\; h(3)= 14,\; h(4)=3.$$

\begin{proposition} \label{prop:cZ} For all $n\ge 5$ one has \quad {\rm (i)} $\Delta C_Z(\Psi_n)=h^{n-4}(3) \quad {\rm (ii)}\; \Delta C_Z(\Lambda_n)=h^{n-5}(3).$
\end{proposition}

\noindent {\it Proof:} The proof is by induction. For $n=5$, we have $\Psi_5= [0,4]$, which has  two points of constancy: $N=1$ and $N=2$.  Therefore $C_Z(\Psi_5)=14=h(3)$. Here the difference 4 is coming from $N=6$, the second point of the interval $\Lambda_5$. We further have $\Lambda_5=[5,7]$, which has  one point of constancy $N=6$. Therefore  $C_Z(\Lambda_5)=3$.

Suppose  the result has been proved till $n$.

\noindent (i) By equation (\ref{eq:rec}),
$$\Delta C_Z(\Psi_{n+1})=  \Delta C_Z(\Psi_n)\Delta C_Z(\Lambda_n)=h^{n-4}(3)h^{n-5}(3)=  h^{n-5}(h(3)3)=h^{n-5}(143)=h^{n-5}(h^2(3))=h^{n-3}(3).$$

\noindent  (ii) Directly from Lemma \ref{lem:cZ}: $\Delta C_Z(\Lambda_{n+1})=\Delta C_Z(\Psi_n)=h^{n-4}(3).$   \hfill\endpf

\bigskip

\noindent {\it Proof of Theorem \ref{th:Zeck-morph}}: The statements on $I_Z$ and $D_Z$ follow immediately from Lemma \ref{lem:diff}.

  The statement on $C_Z$ follows from Proposition \ref{prop:cZ}, part (i), since $h^n(3)=h^n(1)$ for all $n>0$. \hfill\endpf

\section{The  base phi expansion}\label{sec:phi}

 A natural number $N$ is written in base phi (\cite{Bergman}) if $N$ has the form
  $$N= \sum_{i=-\infty}^{\infty} d_i \varphi^i,\vspace*{-.0cm}$$
 with digits $d_i=0$ or 1, and where $d_id_{i+1} = 11$ is not allowed.

 We  write these expansions as
  $$\beta(N) = d_{L}d_{L-1}\dots d_1d_0\cdot d_{-1}d_{-2} \dots d_{R+1}d_R.$$
 Ignoring leading and trailing 0's, the base phi representation of a number $N$ is unique (\cite{Bergman}).

  Let for $N\ge 0$
  $$s_\beta(N):=\sum_{k=L}^{k=R} d_k(N)$$
  be the sum of digits function of the base phi expansions. We have
  $$(s_\beta(N)) =(0, 1, 2, 2, 3, 3, 3, 2, 3, 4, 4, 5, 4, 4, 4, 5, 4, 4, 2, 3, 4, 4, 5, 5, 5, 4, 5, 6, 6, 7, 5, 5, 5, 6,\dots).$$

\medskip

The case of base phi is considerably more complicated than the Zeckendorf case. We need several preparations, before we can prove   Theorem \ref{th:phi-inf} in Section \ref{seq:phi-morph-inf}, Theorem \ref{th:phi-GBS} in Section \ref{sec:phi-GBS} and Theorem \ref{th:phi-morph} in Section \ref{sec:phi-morph}.

\subsection{The Recursive Structure Theorem}\label{sec:luc}

 The  result of this section was anticipated in  \cite{Hart98}, \cite{Hart99}, and \cite{San-San}, and proved in \cite{Dekking-add}.

 The Lucas numbers $(L_n)=(2, 1, 3, 4, 7, 11, 18, 29, 47, 76,123,\dots)$ are defined by
$$  L_0 = 2,\quad L_1 = 1,\quad L_n = L_{n-1} + L_{n-2}\quad {\rm for \:}n\ge 2.$$
For $n\ge 2$ we are interested in  three consecutive intervals given by
\begin{align*}
I_n&:=[L_{2n+1}+1,\, L_{2n+1}+L_{2n-2}-1],\\
J_n&:=[L_{2n+1}+L_{2n-2},\, L_{2n+1}+L_{2n-1}],\\
K_n&:=[L_{2n+1}+L_{2n-1}+1,\, L_{2n+2}-1].
\end{align*}

To formulate the next theorem, it is notationally convenient to extend the semigroup of words to the free group of words. For example, one has $110^{-1}01^{-1}00=100$.

\begin{theorem}{\bf [Recursive Structure Theorem]}\label{th:recstuc}

\noindent{\,\bf I\;} For all $n\ge 1$ and $k=1,\dots,L_{2n-1}$
one has $ \beta(L_{2n}+k) =  \beta(L_{2n})+ \beta(k) = 10\dots0 \,\beta(k)\, 0\dots 01.$

\noindent{\bf II} For all $n\ge 2$ and $k=1,\dots,L_{2n-2}-1$
\begin{align*}
I_n:&\quad \beta(L_{2n+1}+k) = 1000(10)^{-1}\beta(L_{2n-1}+k)(01)^{-1}1001,\\ K_n:&\quad\beta(L_{2n+1}+L_{2n-1}+k)=1010(10)^{-1}\beta(L_{2n-1}+k)(01)^{-1}0001.
\end{align*}
Moreover, for all $n\ge 2$ and $k=0,\dots,L_{2n-3}$
$$\hspace*{0.7cm}J_n:\quad\beta(L_{2n+1}+L_{2n-2}+k) = 10010(10)^{-1}\beta(L_{2n-2}+k)(01)^{-1}001001.$$
\end{theorem}

 It is crucial to our analysis to partition the natural numbers in what we call the Lucas intervals, given by $\Lambda_0:=[0,1]$, and  for $n=1,2\dots$ by
 $$\Lambda_{2n}:=[L_{2n},L_{2n+1}], \quad  \Lambda_{2n+1}:=[L_{2n+1}+1, L_{2n+2}-1]. $$
If $I=[k,\ell]$ and $J=[\ell+1,m]$ are two adjacent intervals of integers, then we write $IJ=[k,m]$.

We code the Lucas intervals with four symbols \bC{0}\,, \bC{1}\,, \bC{2}\, and \bC{3}\,, by a code $\Psi$ in the following way:

\smallskip

\qquad $\Psi(\Lambda_0)=$\bC{0}\,, $\Psi(\Lambda_1)=$\bC{1}\,, $\Psi(\Lambda_2)=$\bC{2}\,, $\Psi(\Lambda_3)=$\bC{3}\,.

\medskip

We then code $\Psi(\Lambda_4)=\Psi(\Lambda_0)\Psi(\Lambda_1)\Psi(\Lambda_2)$=\bC{0}\,\bC{1}\,\bC{2}\,,
 $\Psi(\Lambda_5)=\Psi(\Lambda_3)\Psi(\Lambda_2)\Psi(\Lambda_3)$=\bC{3}\,\bC{2}\,\bC{3}\,,
 and in general  by induction, suggested by Theorem \ref{th:recstuc}:
\begin{align*}
  \Psi(\Lambda_{2n+2}) & =\Psi(\Lambda_0)\Psi(\Lambda_1)\Psi(\Lambda_2)\dots\Psi(\Lambda_{2n}), \\
  \Psi(\Lambda_{2n+1}) & =\Psi(\Lambda_{2n-1})\Psi(\Lambda_{2n-2})\Psi(\Lambda_{2n-1}).
\end{align*}
Let $\sigma$ be the morphism on the alphabet $\{$\bC{0}\,, \bC{1}\,, \bC{2}\,, \bC{3}\,$\}$ defined by

\medskip

\qquad$\sigma$(\bC{0}\,)=\bC{0}\,\bC{1}\,,\quad$\sigma$(\bC{1}\,)=\bC{2}\,\bC{3}\,,
  \quad$\sigma$(\bC{2}\,)=\bC{0}\,\bC{1}\,\bC{2}\,,\quad$\sigma$(\bC{3}\,)=\bC{3}\,\bC{2}\,\bC{3}\,.

\medskip

\begin{lemma}\label{lem:sigma} For each $n\ge 0$ we have \;$\Psi(\Lambda_{2n+2})=\sigma^n($\bC{2}\,$)$,\quad $\Psi(\Lambda_{2n+3})=\sigma^n($\bC{3}\,$)$.
\end{lemma}

\noindent {\it Proof:} By induction. For $n=0$: $\Psi$($\Lambda_{2}$)$=$\bC{2}\,,\;$\Psi(\Lambda_{3})=$\bC{3}\,.
The induction step:

\medskip

 $\Psi$($\Lambda_{2n+5}$) = $\Psi$($\Lambda_{2n+3}$) $\Psi$($\Lambda_{2n+2}$) $\Psi$($\Lambda_{2n+3}$)
= $\sigma^n$(\bC{3}\,)$\sigma^n$(\bC{2}\,)$\sigma^n$(\bC{3}\,) = $\sigma^{n+1}$(\bC{3}\,).

\smallskip

\noindent Also, using the simple identity $\sigma$(\bC{2}\,)\,\bC{3}\,\,$\sigma$(\bC{2}\,) = $\sigma^2$(\bC{2}\,) in the last step:

\smallskip

$\Psi(\Lambda_{2n+4})=\Psi(\Lambda_0)\Psi(\Lambda_1)\Psi(\Lambda_2)\dots\Psi(\Lambda_{2n}) \Psi(\Lambda_{2n+1})\Psi(\Lambda_{2n+2})=
\Psi(\Lambda_{2n+2})\Psi(\Lambda_{2n+1})\Psi(\Lambda_{2n+2})=$

\smallskip

$\sigma^n$(\bC{2}\,)$\sigma^{n-1}$(\bC{3}\,)$\sigma^n$(\bC{2}\,) = $\sigma^{n+1}$(\bC{2}\,) \hfill\endpf

\medskip

We will now show that the fixed point $x_\sigma$ of the morphism $\sigma$ is quasi-Sturmian, and determine its complexity function $p_\sigma$, i.e., $p_\sigma(n)$ is the number of words of length $n$ that occurs in $x_\sigma$. Let $g_{a,b}$ the morphism on the alphabet $\{a,b\}$ given by
\begin{equation}\label{eq:gab}
g_{a,b}(a)=baa, \; g_{a,b}(b)=ba.
\end{equation}
The morphism $g_{a,b}$ is well-known, and closely related to the Fibonacci morphism. In fact, $\xG=bx_{a,b}$, if $\xG$ is the fixed point of $g_{a,b}$, and $x_{a,b}$ is the fixed point of the Fibonacci morphism $a\rightarrow ab, \, b\rightarrow a$ (see \cite{Bers-Seeb}).

\begin{proposition} \label{prop:xs} The fixed point $x_\sigma$ of $\sigma$ is equal to the decoration $\delta(\xG)$ of the fixed point $\xG$ of $g=g_{a,b}$. The decoration morphism $\delta$ is given by $\delta(a)=$\bC{2}\,\bC{3}\,,\, $\delta(b)=$\bC{0}\,\bC{1}\,. \, For all $n\ge 1$ one has $p_\sigma(n)=n+3$.
\end{proposition}

 \noindent {\it Proof:} For the two words \bC{0}\,\bC{1}\, and \bC{2}\,\bC{3}\, occurring in $x_\sigma$ we find
 $$\sigma(\text{\bC{0}\,\bC{1}\,)=\bC{0}\,\bC{1}\,\,\bC{2}\,\bC{3}\,},\quad \sigma(\text{\bC{2}\,\bC{3}\,)=\bC{0}\,\bC{1}\,\,\bC{2}\,\bC{3}\,\,\bC{2}\,\bC{3}\,}.$$
 In other words,
 $$\sigma(\delta(a))=\delta(baa)=\delta(g(a)),\quad \sigma(\delta(b))=\delta(ba)=\delta(g(b)).$$
 Thus   $\sigma\,\delta=\delta \,g$, which implies $\sigma^n\,\delta=\delta \,g^n$ for all $n$. Since $x_\sigma$ has prefix \bC{0}\,\bC{1}\,$\,=\delta(b)$, with $b$ the prefix of $\xG$, this implies the first part of the proposition.

 For the second part, Proposition 8 in \cite{Cassaigne} is not conclusive, as we do not know a priori the constant $n_0$. But there is a direct computation possible. The complexity function of the Sturmian word $\xG$ is given by $p(n)=n+1$. We have, distinguishing between words of even and odd length, and then splitting according to words occurring at even or odd positions in $x_\sigma$,
 $$p_\sigma(2n)=p(n)+p(n+1)=n+1+n+2=2n+3, \quad p_\sigma(2n+1)=p(n+1)+p(n+1)=2n+4.$$
 \hfill\endpf

 \medskip

Proposition \ref{prop:xs} in combination with the main result of  the paper \cite{Huang-Wen}, explains why the factors of $x_\sigma$ have a simple return word structure. This lies at the basis of Theorem \ref{th:phi-GBS} in Section \ref{sec:phi-GBS}.

\subsection{A morphic sequence representation of $s_\beta$}\label{seq:phi-morph-inf}

The image under a morphism $\delta$ of the fixed point $x$ of a morphism, will be called a \emph{decoration} of $x$. It is well known that such a $\delta(x)$ is a morphic sequence, i.e., the letter to letter projection of the fixed point of a morphism. This is the way we formulate the morphic sequence result in the next theorem.

\smallskip

\begin{theorem}\label{th:phi-inf}
  The function $s_\beta$, as a sequence, is a decoration of a morphic sequence on an infinite alphabet, i.e., $(s_\beta(N))$ is an image under a morphism $\delta$ of a fixed point of a morphism $\gamma$.
The alphabet is $\{0,1,...,j,...\}\times\{$\bC{0}$\,,$\bC{1}$\,,$\bC{2}$\,,$\bC{3}$\,\}$, and $\gamma$ is the morphism given for $j\ge 0$ by\\[-.6cm]
\begin{align*}
      \gamma\Big(\PQ{j}{\!$\bC{0}$}\Big)& \;=\; \PQ{j}{\!$\bC{0}$} \PQ{j}{\!$\bC{1}$}, \\
      \gamma\Big(\PQ{j}{\!$\bC{1}$}\Big)& \;=\; \PQ{j}{\!$\bC{2}$} \PQ{j}{\!$\bC{3}$},\\
      \gamma\Big(\PQ{j}{\!$\bC{2}$}\Big)& \;=\; \PQ{j\!+\!2}{$\bC{0}$} \PQ{j\!+\!2}{$\bC{1}$}\PQ{j\!+\!2}{$\bC{2}$},\\
      \gamma\Big(\PQ{j}{\!$\bC{3}$}\Big)& \;=\;  \PQ{j\!+\!1}{$\bC{3}$} \PQ{j\!+\!2}{$\bC{2}$} \PQ{j\!+\!1}{$\bC{3}$}.
\end{align*}
The decoration map is given by the morphism $\delta$:
$$\delta\Big(\PQ{j}{\!$\bC{0}$}\Big)=\gr{0}+j,\gr{1}+j,\quad\delta\Big(\PQ{j}{\!$\bC{1}$}\Big)=\gr{2}+j,
\quad\delta\Big(\PQ{j}{\!$\bC{2}$}\Big)=\gr{2}+j,\gr{3}+j,\quad\delta\Big(\PQ{j}{\!$\bC{3}$}\Big)=\gr{3}+j,\gr{3}+j.  $$
The image $\delta(x_\gamma)$ of the fixed point $x_\gamma$ of $\gamma$ with initial symbol $\pq{0}{\!$\bC{0}$}$ equals $(s_\beta(N))$.
\end{theorem}

\medskip

\noindent {\it Proof:} One combines Theorem \ref{th:recstuc} with Lemma \ref{lem:sigma}. We see from part {\bf I} of Theorem \ref{th:recstuc}, that the number of 1's in the expansion of $N$ from $\Lambda_{2n+2}$ is 2 more than the number of 1's in the corresponding $N'$ in $\Lambda_0\Lambda_1\dots\Lambda_{2n}$. This gives the three upper indices $j+2$ in  $\gamma\Big(\pq{j}{\!$\bC{2}$}\Big)$. Similarly, part {\bf II} gives that the number of 1's in the three intervals
$\Lambda_{2n-1}$, $\Lambda_{2n-2}$, and $\Lambda_{2n-1}$ is increased by 1, by 2, and respectively 1 for the corresponding $N'$ in the interval $\Lambda_{2n+1}$. This gives the three upper indices in  $\gamma\Big(\pq{j}{\!$\bC{3}$}\Big)$. The lower indices are given by the morphism $\sigma$.
This all happens at the level of the shifted versions of the four intervals $\Lambda_0, \Lambda_1,\Lambda_2$ and $\Lambda_3$. Here $\Lambda_0=[0,1]$ with $s_\beta(0)=\gr{0}$ and $s_\beta(1)=\gr{1}$; $\Lambda_1=\{2\}$ with $s_\beta(2)=\gr{2}$; $\Lambda_2=[3,4]$ with $s_\beta(3)=\gr{2}$ and $s_\beta(4)=\gr{3}$; $\Lambda_3=[5,6]$ with $s_\beta(5)=\gr{3}$ and $s_\beta(6)=\gr{3}$. This yields the decorations $\delta$, taking in to account the corresponding increments of the sum of digits.\hfill\endpf

\medskip

We illustrate Theorem \ref{th:phi-inf} with the following table.

\medskip

\hspace*{-.6cm}
\begin{tabular}{|l|cc|c|rl|cc|c|c|}
\hline
 $N^{\phantom{|}}$                          &\!\!\!$0$   & $1$      & $2$      & $3$    &$4$       & $5$ & $6$  & 7\;\;8\;\;9\;\;10\;\,11  &12\;13\;14\;15\;16\;17 \\
\hline
 $s_\beta(N)$                               &\!\!\!$\gr{0}$   & $\gr{1}$ & $\gr{2}$ & \gr{2} & \gr{3}    & \gr{3}    &\gr{3}& \; \gr{2}\;\;\gr{3}\;\;\gr{4}\;\;\;\gr{4}\phantom{0}\;\gr{5}\phantom{1} & \gr{4}\;\;\;\gr{4}\;\;\;\gr{4}\;\;\;\gr{5}\;\;\;\gr{4}\;\;\;\gr{4}\\
\hline
 Lucas interval                             &\phantom{....}$\Lambda_0$&      &$\Lambda_1$&       &$\Lambda_2$\phantom{....}&  &$\Lambda_3$& $\Lambda_4$ & $\Lambda_5$    \\
\hline
 shifted Lucas   intervals                           &\phantom{....}$\Lambda_0$&     &$\Lambda_1$&       &$\Lambda_2$\phantom{....}&  &$\Lambda_3$& $\Lambda_0$\;\;$\Lambda_1$\;\;$\Lambda_2$ & $\Lambda_3$\;\;\;\;$\Lambda_2$\;\;\;\;$\Lambda_3$\\
\hline
\bC{0}\,,\bC{1}\,,\bC{2}\,,\bC{3}\,-coding  &\phantom{....}\bC{0}          &     & \bC{1}    &       &\bC{2}\phantom{....}    &  &\bC{3}     &\bC{0}\;\bC{1}\;\bC{2}\, & \bC{3}\;\bC{2}\;\bC{3}\, \\
\hline
\end{tabular}

\bigskip

\noindent {\bf Remark} In the paper \cite{Dekk-sum-INT} the base phi analogue of the Thue-Morse sequence, i.e., the sequence $(s_\beta(N) \mod 2)$,  is shown to be a morphic sequence. This result follows also from Theorem \ref{th:phi-inf}, by mapping $2j$ to 0, and $2j+1$ to 1. The morphisms found in this way are on a larger alphabet than the morphism in \cite{Dekk-sum-INT}.

\subsection{Generalized Beatty sequences for $s_\beta$}\label{sec:phi-GBS}

Let $I_\beta$ be the sequence listing the points of increase of $s_\beta(N)$ . We see  that the first six points of increase are $I_\beta(1)=0,\; I_\beta(2)=1,\; I_\beta(3)=3,\; I_\beta(4)=7,\; I_\beta(5)=8,\; I_\beta(6)=10$. Similarly we define $C_\beta$ and $D_\beta$.

\medskip

\begin{theorem}\label{th:phi-GBS} The sequence $I_\beta$, the points of increase of the function $s_\beta$,  is given  by  the union of the two generalized Beatty sequences
$$(\lfloor n\varphi \rfloor+2n)_{n\ge 0}, \text{ and }   (4\lfloor n\varphi \rfloor+3n+1)_{n\ge 0}.$$
The sequence $C_\beta$, the points of constancy of the function $s_\beta$,  is given by the union of the four generalized Beatty sequences
$$(3\lfloor n\varphi \rfloor+n+1)_{n\ge 1}, \;(4\lfloor n\varphi \rfloor+3n+2)_{n\ge 0},\;  (7\lfloor n\varphi \rfloor+4n+2)_{n\ge 0}, \text{ and }
(11\lfloor n\varphi \rfloor+7n+4)_{n\ge 1}.$$
The sequence $D_\beta$, the points of decrease of the function $s_\beta$,  is given  by  the union of the three generalized Beatty sequences
$$(4\lfloor n\varphi \rfloor+3n-1)_{n\ge 1}, \;(7\lfloor n\varphi \rfloor+4n)_{n\ge 1}, \text{ and }   (7\lfloor n\varphi \rfloor+4n+4)_{n\ge 1}.$$
\end{theorem}

\medskip

\noindent {\it Proof:} {\bf I: Points of increase}

\smallskip

Any occurrence of a \bC{0}\, gives two points of increase, namely the pair $\gr{0}+j,\gr{1}+j$, and the pair $\gr{1}+j,\gr{2}+j$. Here we use that \bC{0}\, is always followed by \bC{1}\,. Similarly, any occurrence of a \bC{2}\, gives a point of increase $\gr{2}+j,\gr{3}+j$.

As a consequence we obtain the numbers $N$ which are point of increase by the sequences of occurrences of \bC{0}\,, and those of \bC{2}\,. How do we obtain these sequences? We have to study the return words to \bC{0}\,, and \bC{2}\,. The sets of these return words are respectively

\medskip

\qquad $\{$\bC{0}\,\bC{1}\,\bC{2}\,\bC{3}\,,\,  \bC{0}\,\bC{1}\,\bC{2}\,\bC{3}\,\bC{2}\,\bC{3}\,$\}$,
      and $\{$\bC{2}\,\bC{3}\,,\,  \bC{2}\,\bC{3}\,\bC{0}\,\bC{1}\,$\}$.

\medskip

Both \bC{0}\,, and \bC{2}\, induce the descendant morphism $g_{a,b}$.
Here we coded $b$:=\bC{0}\,\bC{1}\,\bC{2}\,\bC{3}\,,\,$a$:=  \bC{0}\,\bC{1}\,\bC{2}\,\bC{3}\,\bC{2}\,\bC{3}\,, respectively
$b$:=\bC{2}\,\bC{3}\,,\,  $a$:=\bC{2}\,\bC{3}\,\bC{0}\,\bC{1}\,.

\smallskip

The occurrences of \bC{0}\, in the fixed point of $\sigma$ occur at distances given by
 the lengths of $\delta($\bC{0}\,\bC{1}\,\bC{2}\,\bC{3}\,) and $\delta($\bC{0}\,\bC{1}\,\bC{2}\,\bC{3}\,\bC{2}\,\bC{3}\,). These are
 $|\delta($\bC{0}\,\bC{1}\,\bC{2}\,\bC{3}\,)$|$=7, and $|\delta($\bC{0}\,\bC{1}\,\bC{2}\,\bC{3}\,\bC{2}\,\bC{3}\,)$|$=11. It then follows from Lemma \ref{lem:diff} that the increase points are given by the union of the two generalized Beatty sequences $V'(4,3,0)$ and $V'(4,3,1)$, where the $'$ indicates that these start from index 0. Similarly, the occurrences of \bC{2} have first differences 7 and 4, giving the  generalized Beatty sequence $V(3,1,-1)$.

 This is not yet the first result in Theorem \ref{th:phi-GBS}, but by Lemma \ref{lem:split} the sequence $V(1,2,0)$ splits into the two sequences $V(3,1,-1)$ and $V(4,3,0)$. Adding $N=0$ to $V(1,2,0)$ and to $V(4,3,0)$ then yields the result on $I_\beta$ in Theorem \ref{th:phi-GBS}.

 \medskip

\noindent   {\bf II: Points of constancy}

\smallskip

 Any occurrence of a \bC{1}\, gives a point of constancy, namely the pair $\gr{2}+j,\gr{2}+j$,  Here we use that \bC{1}\, is always followed by \bC{2}. Similarly, any occurrence of a \bC{3}\, gives a point of constancy $\gr{3}+j,\gr{3}+j$.

 But there are more points of constancy. At the inner boundary of
 $\Lambda_2\Lambda_3$ in the quadruple $\Lambda_0\Lambda_1\Lambda_2\Lambda_3$ occurs $\gr{3},\gr{3}$. However, this is not the case at the inner boundary of
 the interval $\Lambda_2\Lambda_3$ in the triple $\Lambda_3\Lambda_2\Lambda_3$ in $\Lambda_5$.
 Since $\Psi(\Lambda_0\Lambda_2\Lambda_3\Lambda_4)=$\bC{0}\,\bC{1}\,$\sigma($\bC{1}\,$)$, and $\Psi(\Lambda_3\Lambda_2\Lambda_3)=\sigma($\bC{3}\,$)$ these points of constancy occur if and only if $\sigma($\bC{1}\,$)$ occurs in the fixed point of $\sigma$.

 This still does not yet exhaust all possibilities: there is the point $N=14$ with $s_\beta(N)=s_\beta(N+1)=\gr{4}$  in $\Lambda_5$, not yet covered by the previous sequences. This induces points of constancy occurring at all shifted $\Lambda_5$, which occur if and only if $\sigma($\bC{3}\,$)$  occurs in the fixed point of $\sigma$. Since any $\Lambda_k$ for $k>5$ can be written as a union of shifted versions of the three intervals $\Lambda_0\Lambda_1\Lambda_2\Lambda_3$, $\Lambda_4$, and $\Lambda_5$, we have covered all possibilities.

 As a consequence we obtain the numbers $N$ which are point of increase by the sequences of occurrences of \bC{1}\,, \bC{3}\,, $\sigma($\bC{1}\,$)$, and $\sigma($\bC{3}\,$)$. As before, all four have a set of two return words, and a descendant morphism that is equal to $g$.
 For \bC{1}\, the $\delta$-images have lengths 11 and 7, for  \bC{3}\, the $\delta$-images have lengths 7 and 4, for  $\sigma($\bC{1}\,$)$, the $\delta$-images have lengths 29 and 18, and for  $\sigma($\bC{3}\,$)$ the $\delta$-images have lengths 18 and 11. Application of Lemma \ref{lem:diff} then gives the four generalized Beatty sequences of  $C_\beta$ in Theorem \ref{th:phi-GBS}.

 \medskip

{\bf III: Points of decrease}

\smallskip

The first point of decrease is $N=6$, which occurs at the end of $\Lambda_3$,
 so $N+1=7$ occurs at the beginning of $\Lambda_4=\Lambda_0\Lambda_1\Lambda_2$. This gives occurrences of points of decrease at every occurrence of \bC{3}\,\bC{0}\,.  This word has two return words: $b$:=\bC{3}\,\bC{0}\,\bC{1}\,\bC{2}\,, and $a$:=\bC{3}\,\bC{0}\,\bC{1}\,\bC{2}\,\bC{3}\,\bC{2}\,.
 These induce as descendant morphism the morphism $g$, once more. As $|\delta(a)|=11$, and $|\delta(b)|=7$, this leads to the sequence $V'(4,3,-1)$.

 The next point of decrease is at $N=11$, occurring at the inner boundary of the adjacent $\Lambda_4\Lambda_{5}$. The third point of decrease is at $N=15$,
 which lies inside $\Lambda_{5}$. The coding of $\Lambda_{5}$ is $\Psi(\Lambda_{5}) =\, $\bC{3}\,\bC{2}\,\bC{3}\, = $\sigma($\bC{3}\,$)$. As in the previous section, this gives the sequence $V(7,4,0)$ for the occurrences of the decrease points $N=11$, and later shifts. Then $V(7,4,4)$ gives the occurrences of the decrease points $N=15=11+4$, and later shifts. Again, since any $\Lambda_k$ for $k>5$ can be written as a union of intervals $\Lambda_0\Lambda_1\Lambda_2\Lambda_3$, $\Lambda_4$, and $\Lambda_5$, we have covered all possibilities. This finishes the $D_\beta$ part of Theorem \ref{th:phi-GBS}. \hfill\endpf

 \subsection{Morphisms for the first differences }\label{sec:phi-morph}

 As for the Zeckendorf expansion, we have seen in the previous section that the points of constancy have a more complicated structure than the points of increase or the points of decrease. This phenomenon expresses itself also in the `morphic versions' of the characterization.

\begin{theorem}\label{th:phi-morph} The points of increase of the function $s_\beta$  are given by the sequence $I_\beta$, which has $I_\beta(1)=0$, and $\Delta I_\beta$ is the fixed point of the morphism  on the alphabet $\{1,2,4\}$ given by
$$1\rightarrow 12, \; 2\rightarrow 4, \; 4\rightarrow 1244.$$
The points of constancy of the function $s_\beta$  are given by the sequence $C_\beta$, which has $C_\beta(1)=2$, and $\Delta C_\beta$ is a morphic sequence, given by the letter-to-letter projection $1\rightarrow 1, 2\rightarrow 2, 3\rightarrow 3, 3'\rightarrow 3, 4\rightarrow 4$ of the fixed point of the morphism on the alphabet $\{1,2,3,3',4\}$ given by
$$1\rightarrow 43, \; 2\rightarrow 21, \; 3\rightarrow 21,\; 3'\rightarrow 13'43, \; 4\rightarrow 13'4.$$
The points of decrease of the function $s_\beta$  are given by the sequence $D_\beta$, which has $D_\beta(1)=6$, and $\Delta D_\beta$ is the shift by one of the fixed point of the morphism  on the alphabet $\{2,4,5,7\}$ given by
$$2\rightarrow 542,\; 4\rightarrow 542,\;5\rightarrow 7,\;7\rightarrow 7542.$$
\end{theorem}

\noindent {\it Proof:}\,
We use in all three cases  the return words to \bC{0}\, which are $b$:=\bC{0}\,\bC{1}\,\bC{2}\,\bC{3}\, and $a$:=\bC{0}\,\bC{1}\,\bC{2}\,\bC{3}\,\bC{2}\,\bC{3}\, to follow the occurrences of the points of increase, constancy and decrease.
 The important property of these return words is that the \emph{first} occurrence of the points of increase is at the same position in the decorated $a$ and  $b$, and the same holds for the points of constancy and decrease.

\medskip

\noindent {\it Proof:} {\bf I: Points of increase}

\smallskip

We take in to account the increase in the differences of the occurrences of the increase points in  the decorations

\medskip

$\delta\Big(\!\pq{j}{\!$\bC{0}$}\pq{j}{\!$\bC{1}$}\pq{j}{\!$\bC{2}$}\pq{j}{\!$\bC{3}$}\!\Big)=\gr{0}+j,\gr{1}+j,\gr{2}+j,\gr{2}+j,\gr{3}+j,\gr{3}+j,\gr{3}+j,$\\

$\delta\Big(\!\pq{j}{\!$\bC{0}$}\pq{j}{\!$\bC{1}$}\pq{j}{\!$\bC{2}$}\pq{j}{\!$\bC{3}$}\pq{j}{\!$\bC{2}$}\pq{j}{\!$\bC{3}$}\!\Big)=
      \gr{0}+j,\gr{1}+j,\gr{2}+j,\gr{2}+j,\gr{3}+j,\gr{3}+j,\gr{3}+j,\gr{2}+j,\gr{3}+j,\gr{3}+j,\gr{3}+j,$

\smallskip

\noindent of the extended return words $a$ and $b$.  For $a$ these differences  are 1,2,4 and 4. For $b$ the differences between the occurrences of the increase points are 1,2, and 4. Recall here, that the last 4 comes from the first increase point of the next word.  It follows that we can obtain $\Delta I_\beta$ by decorating the fixed point of the morphism $g$ given by $a\rightarrow baa, b\rightarrow ba$ with the two words 124 and 1244. To turn this decorated fixed point into a fixed point, we apply the natural algorithm (cf. the proof of Corollary 9 in \cite{Dekking-deco-TCS}). In this case this gives the following block map on the alphabet $\{a_1, a_2, a_3, a_4,b_1, b_2, b_3\}$:
 \begin{align*}\label{block}
 a_1 a_2a_3a_4 \rightarrow  & \;b_1b_2b_3a_1 a_2a_3a_4a_1 a_2a_3a_4\\
 b_1b_2b_3 \rightarrow      & \;b_1b_2b_3a_1 a_2a_3a_4.
 \end{align*}
 The most efficient way to turn this into a morphism:
 \begin{align*}
 a_1 \rightarrow  &\; b_1b_2, \; a_2 \rightarrow b_3, \; a_3\rightarrow a_1a_2a_3a_4, \; a_4\rightarrow a_1 a_2a_3a_4\\
 b_1 \rightarrow  &\; b_1b_2, \; b_2 \rightarrow b_3, \; b_3 \rightarrow a_1 a_2a_3a_4.
 \end{align*}
 The associated letter-to-letter map $\lambda$ is given by $\lambda(a_1 a_2a_3a_4)=1244,\; \lambda(b_1 b_2b_3)=124.$
 We see that we can consistently merge $a_1$ and $b_1$ to the letter 1, $a_2$ and $b_2$ to the letter 2, and $a_3$ and $b_3$ to the letter 4.
 Renaming $a_4$ by 4, this then yields the morphism $1\rightarrow 12, 2\rightarrow 4, 4\rightarrow 1244$ as generating morphism for $\Delta I_\beta$.

 \medskip

\noindent   {\bf II: Points of constancy}

\smallskip

 We follow the same strategy as in  part {\bf I}. The differences of the occurrences of points of constancy in the decorated versions of $a$ and $b$ are now 2,1,4 and 3,1,3,4. Decorating the fixed point of the morphism $g$ on $\{a,b\}$ by $a\rightarrow 214$, and $b\rightarrow 3134$ this time leads to  a morphism on the alphabet $\{1,2,3,3',4\}$ given by
$$1\rightarrow 43, \; 2\rightarrow 21, \; 3\rightarrow 21,\; 3'\rightarrow 13'43, \; 4\rightarrow 13'4.$$

The letter-to-letter projection $1\rightarrow 1, 2\rightarrow 2, 3\rightarrow 3, 3'\rightarrow 3, 4\rightarrow 4$ of the fixed point of this morphism on the alphabet $\{1,2,3,3',4\}$ yields the sequence $\Delta C_\beta$ (where $C_\beta(1)=2$).

\medskip

\noindent   {\bf III: Points of decrease}

\smallskip

 The differences of the occurrences of points of decrease in the decorated versions of $a$ and $b$ are 7 and 5,4,2. Decorating the fixed point of the morphism $g_{a,b}$  by $a\rightarrow 7$, and $b\rightarrow 542$ this time leads to  a morphism on the alphabet $\{2,4,5,7\}$ given by
$$2\rightarrow 542, \; 4\rightarrow 542, \; 5\rightarrow 7,\; 7\rightarrow 7542.$$

The unique fixed point of this morphism on the alphabet $\{2,4,5,7\}$ yields the sequence $\Delta D_\beta$, when we put $D_\beta(1)=-1$. \hfill\endpf

\section{Alternative proofs of Theorem \ref{th:phi-GBS}  and \ref{th:phi-morph} }\label{sec:phi-alt}


The proofs of Theorem \ref{th:phi-GBS}  and \ref{th:phi-morph} have been based entirely on the properties of the infinite morphism $\gamma$ of Theorem \ref{th:phi-inf}. The question rises whether there is also a more local approach based on the digit blocks of the expansion as was used for the points of constancy, and the points of decrease of the Zeckendorf sum of digits function.
Here we give a sketch of how this might be achieved for the points of increase of the base phi expansion. We say a number $N$ is of type $\B$ if
$d_1d_0d_{-1}(N)=000$, and of type $\E$ if $d_2d_1d_0(N)=001$. One can then prove the following.

\begin{proposition}\label{prop:BE} A number $N$ is a point of increase of $(s_\beta(N))$ if and only if $N$ is of type $\B$ or of type $\E$.\end{proposition}

Next,  Theorem 5.1 from the paper \cite{Dekk-phi-FQ} gives that type $\B$  occurs along the generalized Beatty sequence $(\lfloor n\varphi \rfloor+2n)_{n\ge 0}$, and one can deduce from Remark 6.3 in the same paper that type $\E$ occurs along the generalized Beatty sequence $(4\lfloor n\varphi\rfloor +3n + 1)_{n\ge 0}$. This gives the alternative proof of the $I_B$-part of Theorem \ref{th:phi-GBS}, based on Proposition \ref{prop:BE}.

\bigskip

 We next give a proof of the $\Delta I_{\B}$ part of Theorem \ref{th:phi-morph}, directly from Theorem \ref{th:phi-GBS} by a purely combinatorial argument.

 \medskip

 \noindent {\it Alternative proof of Theorem \ref{th:phi-morph}:}  Let
 $$I_{\B}:=(\lfloor n\varphi\rfloor +2n)_{n\ge 0}, \quad  I_{\E}:=(4\lfloor n\varphi\rfloor +3n + 1)_{n\ge 0}.$$
 By Lemma \ref{lem:diff}, the difference sequence of
 the sequence $(\lfloor n\varphi\rfloor +2n, \,n\ge 1)$ is equal to the Fibonacci word $x_{4,3}=4344344344\dots$ on the alphabet $\{4,3\}$, and the difference sequence of the sequence $(4\lfloor n\varphi\rfloor +3n + 1,\,n\ge 1)$ is the Fibonacci word $x_{11,7}=11,7,11,11,7,\dots$.
 However, in Theorem \ref{th:phi-GBS} the sequences start at $n=0$, yielding the two difference sequences
 $$\Delta I_{\B}=3x_{4,3}=34344344344\dots, \quad \Delta I_{\E}=7x_{11,7}=7,11,7,11,11,7,\dots.$$
  Recall that the sequences $bx_{a,b}$ are fixed points of the morphisms $g_{a,b}$ from Equation (\ref{eq:gab}) given by
 $g_{a,b}(a)=baa, \;g_{a,b}(b)=ba .$
 The return words of 3 in $\Delta I_{\B}$ are 34 and 344. We code these words by the differences that they yield between successive occurrences of 3's, i.e., by the letters 7 and 11.   Then, since
 $$g_{4,3}(34) =34\,344, \quad g_{4,3}(344) =34\,344\, 344,$$
  the return words induce a derived morphism
 $$7\rightarrow 7,11,\quad 11\rightarrow 7,11,11.$$
 This derived morphism happens to be equal to $g_{11,7}$, the morphism giving the sequence $\Delta I_{\E}$. This implies that to merge the two sequences $I_{\B}$ and $I_{\E}$ to obtain $I$, one has to replace the 3's in $\Delta I_{\B}$ by 1,2. This  decoration of $\Delta I_{\B}$, induces a morphism $\mu$ on the alphabet $\{1,2,4\}$ in the usual way, given by
 $$ \mu(1)= 12,\;\mu(2)=4,\; \mu(4)=1244.$$
 This proves the theorem. \hfill\endpf


\begin{thebibliography}{99}




\bibitem{GBS} J.-P.~Allouche and F.M.~Dekking, Generalized Beatty sequences and complementary triples,
Mosc.~J.~Comb.~Number Theory 8 (2019), 325--342. \; doi.org/10.2140/moscow.2019.8.325

\bibitem{Bergman}
G.~Bergman, A number system with an irrational base, Math. Mag. 31 (1957), 98--110.

\bibitem{Bers-Seeb} J. Berstel, P. S$\acute{e}\acute{e}$bold, A remark on morphic Sturmian words, RAIRO Theor. Inform. Appl. 28 (1994) 255--263.

\bibitem{Cassaigne} J. Cassaigne, Sequences with grouped factors, in {\it DLT'97,
Developments in Language Theory {III}, Thessaloniki, Aristotle University of Thessaloniki}, 1998,
pp.~211--222.

		

%
%
%
%
%
%

\bibitem{Dekking-deco-TCS} F.M. Dekking,  Morphic words, Beatty sequences and integer images of the Fibonacci language, Theoretical Computer Science   809 (2019), 407--417.   \;doi.org/10.1016/j.tcs.2019.12.036


\bibitem{Dekk-phi-FQ} F.M. Dekking, Base phi representations and golden mean beta-expansions, Fibonacci Quart. 58  (2020), 38--48.

\bibitem{Dekk-sum-INT} M. Dekking, The sum of digits function of the base phi expansion of the natural numbers, INTEGERS \underline{20},  \#A45  (2020), 1--6.



\bibitem{Dekking-add} F. Michel Dekking, How to add two natural numbers in base phi. To appear in  Fibonacci Quart. (2020).

 \bibitem{Dekking-Zeck} F. Michel Dekking, The structure of Zeckendorf expansions. arXiv:2006.06970v1  (2020).



%
%
%
%
%
%
%
%
%
%
%
%
%


\bibitem{Hart98} E.~Hart, On Using Patterns in the Beta-Expansions To Study Fibonacci-Lucas Products,
 Fibonacci Quart. 36 (1998), 396--406.


\bibitem{Hart99} E.~Hart  and  L.~Sanchis, On the occurrence of $F_n$ in the Zeckendorf decomposition of $nF_n$,
Fibonacci Quart. 37 (1999), 21--33.

\bibitem{Holton-Zamboni} C.~Holton and L.~Q.~Zamboni, Descendants of primitive substitutions, Theory Comput. Systems 32 (1999), 133–-157.

\bibitem{Huang-Wen} Y.~Huang, Z.Y.~Wen, The sequence of return words of the Fibonacci sequence, Theoretical Computer Science
 593,  106--116.



\bibitem{lothaire} M.~Lothaire, {\it Algebraic combinatorics on words}, Cambridge
University Press, 2002.






\bibitem{oeis} {\it On-Line Encyclopedia of Integer Sequences}, founded by N. J. A. Sloane,
electronically available at {http://oeis.org}.




\bibitem{San-San}  G.R.~Sanchis and L.A.~Sanchis,
On the frequency of occurrence of $\alpha^i$ in the $\alpha$-expansions of the positive integers, Fibonacci Quart. 39 (2001), 123--173.









\end{thebibliography}
\end{document}